\documentclass{article}

\usepackage{mathrsfs}
\usepackage{amsmath, amsthm, amssymb, amsfonts}
\usepackage{graphicx}

\usepackage{texdraw}    
\def\Rad{0.1}
\def\factor{0.6}
\def\unit{cm}
\newdimen\Radius
\Radius=\Rad\unit
\Radius=\factor\Radius
\everytexdraw{%
\drawdim pt
\linewd 0.2
\drawdim{\unit}
\setunitscale{\factor}
}
{\catcode`\p=12\catcode`\t=12\gdef\sri#1pt{#1}}
\newbox\txt
\newdimen\Hoehe
\def\Zeichen#1{\htext{\BOX{#1}}}
\def\BOX#1{\setbox\txt=\vbox{\baselineskip6pt%
\halign{$\scriptstyle##$\hfil\cr#1\cr}}%
\global\edef\B{\expandafter\sri\the\wd\txt}%
\global\Hoehe=\ht\txt
\global\advance\Hoehe\dp\txt
\global\edef\H{\expandafter\sri\the\Hoehe}%
\copy\txt}
\def\punkt(#1,#2){\move(#1 #2)\fcir f:0 r:\Rad}
\def\kreis(#1,#2){\move(#1 #2)\lcir r:\RRad}
\def\point(#1,#2){\move(#1 #2)\fcir f:0 r:\Rad}
\def\circle(#1,#2){\move(#1 #2)\lcir r:\RRad}
\def\strichRO(#1,#2){\move(#1 #2)\rlvec(1 1)}
\def\strichRU(#1,#2){\move(#1 #2)\rlvec(1 -1)}
\def\strichU(#1,#2){\move(#1 #2)\rlvec(0 -1)}
\def\lined(#1,#2){\move(#1 #2)\rlvec(0 -1)}
\def\linedl(#1,#2){\move(#1 #2)\rlvec(-1 -1)}
\def\linedr(#1,#2){\move(#1 #2)\rlvec(1 -1)}
\def\linedll(#1,#2){\move(#1 #2)\rlvec(-2 -1)}
\def\linedrr(#1,#2){\move(#1 #2)\rlvec(2 -1)}
\def\linedlll(#1,#2){\move(#1 #2)\rlvec(-3 -1)}
\def\linedrrr(#1,#2){\move(#1 #2)\rlvec(3 -1)}
\def\labelO(#1,#2)[#3]#4{%
\punkt(#1,#2)
\rmove(#3 0)
\textref h:C v:B
\Zeichen{#4\vrule depth 5pt width 0pt}
{\drawdim pt
\setunitscale 1
\rmove(0 {\H})
\setunitscale 0.5
\rmove({\B} 0)
\setunitscale 1
\rmove({-\B} 0)
}}
\def\labelU(#1,#2)[#3]#4{%
\punkt(#1,#2)
\rmove(#3 -0.15)
\textref h:C v:T
\Zeichen{#4}
{\drawdim pt
\setunitscale 1
\rmove(0 {-\H})
\setunitscale 0.5
\rmove({\B} 0)
\setunitscale 1
\rmove({-\B} 0)
}\ifx0#2\realmult{0.15}{\factor}\fff
\advance\Hoehe by \fff\unit
\ifdim\Hoehe > \maxT \global\maxT=\Hoehe\fi\fi}
\def\labelR(#1,#2)[#3]#4{%
\punkt(#1,#2)
\rmove(.2 #3)
\textref h:L v:C
\Zeichen{#4}
{\drawdim pt
\setunitscale 1
\rmove({\B} 0)
\setunitscale 0.5
\rmove(0 {\H})
\setunitscale 1
\rmove(0 {-\H})
}\ifx0#2\Hoehe=0.5\Hoehe
\realmult{#3}{\factor}\fff
\advance\Hoehe by -\fff\unit
\ifdim\Hoehe > \maxT \global\maxT=\Hoehe\fi\fi}
\def\labelL(#1,#2)[#3]#4{%
\punkt(#1,#2)
\rmove(-.2 #3)
\textref h:R v:C
\Zeichen{#4}
{\drawdim pt
\setunitscale 1
\rmove({-\B} 0)
\setunitscale 0.5
\rmove(0 {\H})
\setunitscale 1
\rmove(0 {-\H})
}\ifx0#2\Hoehe=0.5\Hoehe
\realmult{#3}{\factor}\fff
\advance\Hoehe by -\fff\unit
\ifdim\Hoehe > \maxT \global\maxT=\Hoehe\fi\fi}
\def\labelcond(#1,#2)[#3]#4{%
\move(#1 #2)
\rmove(#3 -0.15)
\textref h:C v:T
\Zeichen{#4}
{\drawdim pt
\setunitscale 1
\rmove(0 {-\H})
\setunitscale 0.5
\rmove({\B} 0)
\setunitscale 1
\rmove({-\B} 0)
}\ifx0#2\realmult{0.15}{\factor}\fff
\advance\Hoehe by \fff\unit
\ifdim\Hoehe > \maxT \global\maxT=\Hoehe\fi\fi}
\def\Vstern(#1,#2){%
\move(#1 #2)
\rmove(-0.1 0.5)
\textref h:R v:C
\Zeichen{\textstyle*}
{\drawdim pt
\setunitscale 1
\rmove({-\B} 0)
\setunitscale 0.5
\rmove(0 {\H})
\setunitscale 1
\rmove(0 {-\H})
}}
\def\Lstern(#1,#2){%
\move(#1 #2)
\rmove(-0.8 0.4)
\textref h:C v:C
\Zeichen{\textstyle*}}
\def\Mstern(#1,#2){%
\move(#1 #2)
\rmove(0 -0.4)
\textref h:C v:C
\Zeichen{\textstyle*}}
\newdimen\maxT
\newbox\Zbox
\long\def\bdiadraw#1\ediadraw{%
\maxT=0pt\setbox\Zbox=\hbox{\btexdraw #1\etexdraw}%
\ifdim\maxT>\Radius\advance\maxT by -\Radius\fi
\lower\maxT\box\Zbox}
\def\labelOO(#1,#2)[#3]#4{%
\move(#1 #2)
\rmove(#3 0)
\textref h:C v:B
\Zeichen{#4\vrule depth 5pt width 0pt}
{\drawdim pt
\setunitscale 1
\rmove(0 {\H})
\setunitscale 0.5
\rmove({\B} 0)
\setunitscale 1
\rmove({-\B} 0)
}}
\def\labelUU(#1,#2)[#3]#4{%
\move(#1 #2)
\rmove(#3 -0.15)
\textref h:C v:T
\Zeichen{#4}
{\drawdim pt
\setunitscale 1
\rmove(0 {-\H})
\setunitscale 0.5
\rmove({\B} 0)
\setunitscale 1
\rmove({-\B} 0)
}\ifx0#2\realmult{0.15}{\factor}\fff
\advance\Hoehe by \fff\unit
\ifdim\Hoehe > \maxT \global\maxT=\Hoehe\fi\fi}
\def\labelUu(#1,#2)[#3]#4{%
\punkt(#1,#2)
\rmove(#3 -0.15)
\textref h:C v:T
\Zeichen{#4}
{\drawdim pt
\setunitscale 1
\rmove(0 {\H})
\setunitscale 0.5
\rmove({\B} 0)
\setunitscale 1
\rmove({-\B} 0)
}}
\def\labelRR(#1,#2)[#3]#4{%
\move(#1 #2)
\rmove(.2 #3)
\textref h:L v:C
\Zeichen{#4}}
\def\labelLL(#1,#2)[#3]#4{%
\move(#1 #2)
\rmove(-.2 #3)
\textref h:R v:C
\Zeichen{#4}}

\renewcommand{\P}{\mathbf P}           
             
\newcommand{\Grass}{\operatorname{G}}      
\newcommand{\poset}{{\mathcal P}}      
\newcommand{\N}{{\mathcal N}}          
\newcommand{\HH}{\mathrm H}            
\newcommand{\I}{{\mathcal I}}          
\renewcommand{\O}{{\mathcal O}}        
\newcommand{\Hom}{\mathrm {Hom}}       
\newcommand{\sheafHom}{{\mathcal H}om} 
\newcommand{\x}{\mathbf{x}}            
\newcommand{\y}{\mathbf{y}}            
\renewcommand{\a}{\alpha}

\renewcommand{\l}{\lambda}

\newcommand{\M}{{\mathcal M}}
\newcommand{\du}{\stackrel{.}{\cup}}
\newcommand{\e}{\varepsilon}
\newcommand{\Z}{\mathbf{Z}}            
\newcommand{\F}{\mathcal{F}}           
\newcommand{\R}{\mathcal{R}}           
\newcommand{\G}{\mathcal{G}}         
\newcommand{\reg}{\mathrm{reg}}

\newcommand{\Hilb}{{\mathscr{H}}}       
\newcommand{\HP}{{\operatorname{Hilb}}} 
\newcommand{\poly}{\mathbf{p}}          
\newcommand{\D}{\Delta}
\newcommand{\supp}{\mathrm{supp}}       
\newcommand{\init}{\operatorname{in}}   
\newcommand{\sat}{\operatorname{sat}}   
\newcommand{\GL}{{\operatorname{GL}}}   
\newcommand{\geb}{{\geq_{\mathrm{B}}}}  
\newcommand{\leb}{{\leq_{\mathrm{B}}}}  

\newtheorem{theorem}{Theorem}[section]
\newtheorem{proposition}[theorem]{Proposition}

\newtheorem{lemma}[theorem]{Lemma}
\newtheorem{corollary}[theorem]{Corollary}

\begin{document}

\begin{center}
\Large \bf On an extension of Galligo's theorem concerning the 
		Borel-fixed points on the Hilbert scheme.
\end{center}

\begin{center}
Morgan Sherman\\
California State University, Channel Islands
\end{center}

\begin{center}
\begin{minipage}{4in}
\centerline{\bf Abstract}

\small
Given an ideal $I$ and a weight vector $w$ which partially orders 
monomials we can consider the initial ideal $\init_w (I)$ which 
has the same Hilbert function.  A well known construction 
carries this out via a one-parameter subgroup of a 
$\GL_{n+1}$ which can then be viewed as a curve on the 
corresponding Hilbert scheme.  Galligo \cite{galligo} 
proved that if $I$ is in generic coordinates, and if 
$w$ induces a monomial order up to a large enough degree, 
then $\init_w(I)$ is fixed by the action of the Borel 
subgroup of upper-triangular matrices.  We prove that the 
direction the path approaches this Borel-fixed point 
on the Hilbert scheme is also Borel-fixed.  
\end{minipage}
\end{center}


\section{Introduction}

The purpose of this paper is to prove a first order 
infinitesimal version of a theorem of Galligo 
\cite{galligo}.  Galligo's theorem states that 
in generic coordinates the initial ideal of any ideal 
is fixed by the action of the Borel 
subgroup of invertible upper-triangular matrices.  This 
theorem has important consequences when translated 
to the Hilbert scheme.  For example it immediately 
follows that any component, and any intersection of 
components on the Hilbert scheme will contain a 
Borel-fixed {\it point}.  This follows since once 
we associate an ideal with its corresponding point on the 
Hilbert scheme, taking an initial ideal corresponds 
(in a way we make precise below) to the closure of 
a path paremetrized 
by an appropriate one-parameter subgroup of a $\GL_{n+1}.$  

The infinitesimal version proven here says that not only is 
the limit point of this path Borel-fixed (Galligo's theorem 
translated to the Hilbert scheme), but also the path 
picks out a vector in the tangent space of the limit point 
which spans a subspace which 
is itself Borel-fixed.  (Note that since the limit 
point is Borel-fixed, the action of the Borel group will 
descend to an action on the tangent space.)  

This problem was posed to me by my PhD advisor 
David Bayer at Columbia University.  I am greatful to him 
for many helpful conversations.  The problem is also part 
of an on-going project to understand the local structure of 
the Hilbert scheme at a Borel-fixed point.  

This paper is organized as follows:  In sections 
\ref{The Hilbert scheme} 
and 
\ref{Borel-fixed ideals} 
we quickly reproduce the relevant information 
needed about Hilbert schemes and Borel-fixed ideals.  
In section 
\ref{The poset P(m,n)} 
we introduce a poset designed to 
capture combinatorially all the information of a 
Borel-fixed ideal.  We discuss the poset and some of 
its properties briefly.  The author believes the poset 
is in some ways the proper way to think about Borel-fixed 
ideals.  Indeed, using the language of posets significantly 
eases statements of the later theorems.  In section 
\ref{Description of the tangent space} 
we develop the notation used for the tangent space to the 
Hilbert scheme at a Borel-fixed point.  Sections 
\ref{The maximal torus eigenvectors} 
and 
\ref{The Borel eigenvectors} 
classify all the vectors of the tangent space which are 
Borel-eigenvectors (span a Borel-fixed subspace).  
Finally, in section 7 we prove the main result of this paper.  


\section{The Hilbert scheme}\label{The Hilbert scheme}

Throughout this paper we 
will work over an algebraically closed field $K$ of 
characteristic $0.$  
Let $\Hilb^{\P^n}_{\poly(z)}$ (or simply $\Hilb$) 
denote the Hilbert Scheme parametrizing all 
subschemes of $\P^n$ with a fixed Hilbert Polynomial
$\poly(z)$.  We set $S = K [ x_0, \ldots, x_n ]$ to be the 
homogeneous coordinate ring for $\P^n$, and for $d \geq 0$ 
we denote by 
$S_d$ the vector space of the homogeneous forms of 
degree $d$ in $S$, so that $S = \oplus_{d \geq 0} S_d$.  
Similarly
for any homogeneous ideal $I \subseteq S$ we denote by
$I_d$ the vector space of its $d$th graded piece.  
Furthermore $I_{\geq d}$ denotes the truncated ideal with 
all elements of degree less than $d$ removed.  
If $f_1, \ldots, f_r \in S$ we will write 
$ ( f_1, \ldots, f_r ) $ for the ideal generated by the
$f_i$s. 

The group $\operatorname{GL}(n+1, K)$ acts on $S$ by 
extending its action on $S_1 \cong K^{n+1}.$  The 
action on $S_1$ is computed in matrix form by taking 
$\{x_0, \ldots, x_n\}$ to be a basis.  
If $g = (a_{ij})$ then 
\[ 
   g(x_i) = g \cdot x_i = 
   a_{0i} x_0 + \ldots + a_{ni} x_n.  
\]
For a simple example, if 
$ g = \left(\begin{array}{cc}1&1\\ 0&1\end{array}\right)$
and $S = k[x,y]$ (where we set $x = x_0, y = x_1$) 
then $g(x) = x$ while $g(y) = x+y$ and hence $g(xy) = x^2 + xy.$  
$\operatorname{GL}(n+1, K)$ then acts on the set of ideals 
of $S$ as well.  Furthermore, if $I$ is an ideal of $S$ then 
$g(I)$ defines a scheme projectively equivalent to that 
defined by $I,$ so we get an action of $\operatorname{GL}(n+1, K)$ 
on $\Hilb,$ and in fact on each of its irreducible components. 

Given a scheme $Z \subseteq \P^n$, there are many ideals which define
it.  Among all such ideals there is a unique maximal one which 
contains all others.  It can be obtained by the 
global sections functor $\I \mapsto
\oplus_{d \geq 0} \HH^0(\I (d))$ applied 
to the ideal sheaf $\I$ of $Z,$ or 
equivalently by taking the primary decomposition of 
any ideal defining $Z$ and removing the component 
associated to $(x_0,\ldots,x_n).$  This operation is called 
\emph{saturation} and the result of applying it to the ideal 
$I$ will be denoted as $I^{\sat}.$  
 
Any two ideals defining $Z$ will agree in large enough 
degree.  Thus if $I$ is the saturated ideal defining $Z$ and 
$d$ is large enough then the $d$-th 
graded piece $I_d$ determines $Z$:  the ideal generated 
by $I_d$ agrees with $I$ in degrees $d$ and above;  saturating 
the result recovers $I.$  
 
The question of how large $d$ should be is answered in part by noting 
the \emph{regularity} of $I$ suffices.  We briefly 
recall what this is (the notion of regularity is due to 
Castelnuovo and Mumford;  for a more detailed account see 
\cite{mumford}).  For any coherent sheaf 
$\mathscr{F}$ and nonnegative integer $m,$ we 
say $\mathscr{F}$ is $m-$regular if 
$\HH^i \mathscr{F}(m-i) = 0$ for all $i > 0.$  
The regularity of $\mathscr{F}$ is the least 
integer $m$ for which $\mathscr{F}$ is $m-$regular.  
Castelnuovo proved that if $\mathscr{F}$ is $m-$regular, then: 
(i)  $\mathscr{F}$ is $j-$regular for each $j \geq m;$  
(ii) $\mathscr{F}(m)$ is generated by global sections.  
The regularity can also be characterized in terms of a 
minimal free resolution.  This has the benefit of 
allowing one to define regularity for any finitely generated 
module.  Let 
\[ 
   0 \leftarrow F \leftarrow \bigoplus_j S(-e_{0j}) \leftarrow 
   \cdots \leftarrow \bigoplus_j S(-e_{nj}) \leftarrow 0 
\]
be a minimal graded free resolution of the finitely generated module 
$F.$  Then the regularity of $F$ is 
$\max \{e_{ij}-i\}.$  

It is very convenient that there is a finite integer bounding the 
regularity of all saturated ideals of schemes 
with a fixed Hilbert polynomial  \cite{gotzmann}.  
The smallest such integer is known as the \emph{Gotzmann number}.  

Given a Hilbert polynomial $\poly(z)$ the Gotzmann number can be 
readily computed.  Write $\poly(z)$ in the form 
\begin{equation}\label{eqn: g polynomial}
   \poly(z) = g(m_0, \ldots, m_s; z) := \sum_{i=0}^{s} 
   {z+i \choose i+1} - {z+i-m_i \choose i+1} . 
\end{equation}
(See \cite{macaulay} for details.)
The integers $m_0, \ldots, m_s$ satisfy 
$m_0 \geq m_1 \geq \ldots \geq m_s$ and 
are unique if we require $m_s \neq 0,$ in which 
case we also get $s \leq n$ (in fact $s$ is the 
dimension of the scheme).    
The Gotzmann number can be read off as $m_0.$
 
Fix a Hilbert polynomial $\poly(z)$ and let $m$ be the Gotzmann 
number.  Since a scheme with Hilbert polynomial $\poly(z)$ can 
be identified with the vector space of degree $m$ forms in 
its saturated defining ideal we can make a set-theoretical 
identification 
\[ 
   \Hilb \cong \{ I_m \mid I = I^{\sat}, \ 
   \poly_{S/I}(z) = \poly(z) \}. 
\]
Let $s = \dim S_m$ and $r = s - \poly(m).$  Any $I_m$ in the 
above set has dimension equal to $r,$ and is a subspace 
of $S_m.$  This gives a set-theoretical inclusion of the above 
set into $\operatorname{G}(r, S_m),$ the Grassmanian of $r$-dimensional 
subspaces of $S_m.$  Thus we have a set-thoeretical inclusion 
of $\Hilb$ into $\operatorname{G}(r, S_m).$ 
One only needs to verify that this inclusion 
identifies $\Hilb$ with a closed subscheme of 
$\operatorname{G}(r, S_m)$ with the proper scheme structure.  
This is accomplished by using the equations arising from 
the condition
\begin{gather*}
	V \in \{ I_m \mid I = I^{\sat}, \ \HP(S/I) = \poly(z) \} \\
	\Updownarrow   \\
	\dim \left( \{\mathrm{Ideal\ generated\ by\ V}\}_{m+1} \right)
		= \dim S_{m+1} - \poly(m+1).    
\end{gather*}
The fact that these equations scheme-theoretically 
define the Hilbert scheme was conjectured by 
Bayer \cite{bayer} and proven by 
Haiman and Sturmfels \cite{haiman/sturmfels}.  

Throughout the paper we will be viewing the Hilbert 
scheme in this way;  the points will correspond with 
vector subspaces of the vector space $S_m,$ for an 
appropriately chosen $m.$


\section{Borel-fixed ideals}\label{Borel-fixed ideals}

If $A = (a_0, \ldots, a_n) \in \mathbf{N}^{n+1}$ is a vector 
of non-negative integers, and $\x = (x_0, \ldots, x_n)$, then we 
use the notation $\x^A$ for the monomial 
$x_0^{a_0} \cdots x_n^{a_n}$.  We will refer to $A$ as 
the exponent vector of $\x^A.$  

A \emph{monomial order} (or \emph{term order}) is 
a total multiplicative order on the set of monomials 
such that $1$ is the least monomial.  If 
$S = K[x_0,\ldots,x_n]$ we will assume throughout that 
any monomomial order $>$ satisfies 
$x_0 > x_1 > \cdots > x_n.$  

If $I$ is a monomial ideal of $S = K[x_0,\ldots,x_n]$ 
(that is the minimal non-zero generators of $I$ are 
monomials;  equivalently $I$ is fixed by the action 
diagonal matrices in $\GL(n+1,K)$), we set 
$\M(I)$ to be the set of monomials lying in 
$I,$ and $\M(I_d)$ the set of monomials lying in 
$I_d$ (the degree $d$ monomials of $I$).  
Also $\G(I)$ will denote
the minimal generating set of monomials for $I.$ 
For a monomial $\x^A$ we set 
$\max (\x^A )$ (or simply $\max(A)$) to be the index of the 
last variable dividing $\x^A.$  
That is 
\[ \max (\x^A ) = \max(A) := \max \{i \mid x_i | \x^A \} \]
We similarly define $\min(\x^A)$ (and $\min(A)$).  Note that 
with this definition it makes sense to set 
$\max(1) = -\infty,$ and $\min(1) = +\infty.$ 
Finally, $\deg(x^A)$ (or $\deg(A)$) denotes the degree 
of the monomial $\x^A,$ and 
$\deg_i(\x^A)$ (or $\deg_i(A)$) denotes the degree to 
which the variable $x_i$ appears in $\x^A.$  

Recall the action of $\GL(n+1,K)$ on the set of ideals of 
$S = K[x_0,\ldots,x_n].$  An ideal is said to be 
\emph{Borel-fixed} if it is fixed by the action of the 
Borel subgroup of $\GL(n+1,K)$ consisting of upper triangular 
matrices.  Such ideals are stable ideals (defined in the next section)  
in the sense of Eliahou and Kervaire \cite{eliahou/kervaire} and 
are precisely the strongly stable ideals 
in the sense of Peeva and Stillman 
\cite{peeva/stillman}.  Their 
corresponding points on the Hilbert scheme are of significant geometrical 
importance by virtue of their fixed point status.  
Moreover these ideals can be easily classified.  

\begin{proposition}\label{Borel-fixed ideal classification}
The Borel-fixed ideals are the ideals
$I$ such that \\
(1)  $I$ is a monomial ideal. \\
(2)  If $\x^A \in I$ is a monomial, and $x_j \mid \x^A,$
then for $i < j, \ \frac{x_i}{x_j}\x^A \in I$  
\end{proposition}

See, for example, \cite{eisenbud} chapter 15.  
The saturation and the regularity of a 
Borel-fixed ideal are easy to 
determine:

\begin{theorem}
Let $I$ be a Borel-fixed ideal, and 
$\G(I)$ its set of minimal monomial generators. 
Then the saturation of the ideal 
is generated by $\G(I)|_{x_n = 1},$  that is 
one deletes the variable $x_n$ in each of the generators.
\end{theorem}

\begin{proof}
A monomial $\x^A$ is in $I^{\sat}$ iff there is a power 
$k$ such that $x_i^k \x^A \in I$ for $i = 0, \ldots, n.$  
Since $I$ is Borel-fixed this is the case iff 
$x_n^k \x^A \in I$ 
(proposition \ref{Borel-fixed ideal classification}).  
Thus for any monomial $\x^A \in S,$ we 
find $\x^A \in I^{\sat}$ if and only if there is some 
monomial $\x^B \in \G(I)$ 
such that $\x^B$ divides $\x^A x_n^k$ for $k \gg 0.$  One 
sees that this is equivalent to $\x^A$ being a multiple 
of $\x^B\!\mid_{x_n=1}.$  
\end{proof}

\begin{theorem}\label{thm: regularity}
The regularity of a 
Borel-fixed ideal is the highest degree 
of its minimal monomial generators.
\end{theorem}

\begin{proof}  See \cite{bayer}.  \end{proof}

For any ideal $I$ and term order $>,$ the 
\emph{initial ideal} $\init_> I$ is the 
monomial ideal generated by the largest monomials 
appearing in all polynomials of $I.$  The initial 
ideal can be obtained as a $1-$parameter flat 
deformation of $I$ (see section  
\ref{An infinitesimal version of Galligo's theorem}).  
A theorem of Bayer and Stillman and \cite{bayer/stillman} 
states that \emph{in generic coordinates, 
the regularity of an ideal is 
equal to the regularity of its initial ideal in the 
reverse lexicographic order}.  
Thus theorem 
\ref{thm: regularity} takes on great significance in 
the problem of determining regularity.  


\section{The poset $\poset(m,n)$}\label{The poset P(m,n)}

Proposition \ref{Borel-fixed ideal classification} 
endows a Borel-fixed ideal with a combinatorial structure.  
Let $\poset = \poset(m,n)$ be the poset on the set 
of monomials of degree $m$ in $S = k[x_0,\ldots,x_n]$ with the relation 
$\geb$ generated by the covering relation 
$\succ_{\mathrm{B}}$ where
\[ 
   \x^A \succ_{\mathrm{B}} \x^B \iff 
   \exists i < n \mathrm{\ such\ that\ } 
   \x^A = \frac{x_i}{x_{i+1}} \x^B  
\]
We note that \emph{every monomial order $>$ satisfying 
$x_0 > x_1 > \ldots > x_n$ 
is a refinement of this Borel (partial) order.} 
Similar posets are considered in \cite{marinari/ramella} and 
\cite{snellman}.  

For any Borel-fixed ideal $I,$ the set 
$\M(I_m)$ of monomials in $I_m$ will
constitute a {\it filter} of $\poset(m,n)$ -- that is a subset  
$\F \subseteq \poset$ such that 
$\x^B \in \F$ and $\x^A \geb \x^B$ implies $\x^A \in \F$.
Dually, the {\it standard monomials} of degree $m$ for
$I$ (monomials in $S_m \setminus I_m$) constitute
an {\it order ideal} of $\poset$, that is a subset
$\R \subseteq \poset$ such that if
$\x^B \in \R$ and $\x^A \leb \x^B$ then $\x^A \in \R$.  

\begin{figure}[h]
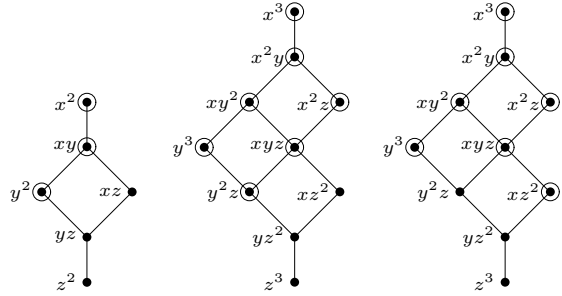

\begin{center}\begin{texdraw}
\labelL(1,4)[0]{x^2} \strichU(1,4) \kreis(1,4)
\labelL(1,3)[0]{xy}  \strichRU(1,3)\kreis(1,3)
\labelL(2,2)[0]{xz}  
\labelL(0,2)[0]{y^2} \strichRU(0,2)\strichRO(0,2)\kreis(0,2)
\labelL(1,1)[0]{yz}  \strichU(1,1)\strichRO(1,1)
\labelL(1,0)[0]{z^2}
\end{texdraw}
\quad
\begin{texdraw}
\labelL(0,7)[0]{x^3}\lined(0,7)\circle(0,7)
\labelL(0,6)[0]{x^2y}\linedl(0,6)\linedr(0,6)\circle(0,6)
\labelL(-1,5)[0]{xy^2}\linedl(-1,5)\linedr(-1,5)\circle(-1,5)
\labelL(1,5)[0]{x^2z}\linedl(1,5)\circle(1,5)
\labelL(-2,4)[0]{y^3}\linedr(-2,4)\circle(-2,4)
\labelL(0,4)[0]{xyz}\linedl(0,4)\linedr(0,4)\circle(0,4)
\labelL(-1,3)[0]{y^2z}\linedr(-1,3)\circle(-1,3)
\labelL(1,3)[0]{xz^2}\linedl(1,3)
\labelL(0,2)[0]{yz^2}\lined(0,2)
\labelL(0,1)[0]{z^3}
\end{texdraw}
\quad
\begin{texdraw}
\labelL(0,7)[0]{x^3}\lined(0,7)\circle(0,7)
\labelL(0,6)[0]{x^2y}\linedl(0,6)\linedr(0,6)\circle(0,6)
\labelL(-1,5)[0]{xy^2}\linedl(-1,5)\linedr(-1,5)\circle(-1,5)
\labelL(1,5)[0]{x^2z}\linedl(1,5)\circle(1,5)
\labelL(-2,4)[0]{y^3}\linedr(-2,4)\circle(-2,4)
\labelL(0,4)[0]{xyz}\linedl(0,4)\linedr(0,4)\circle(0,4)
\labelL(-1,3)[0]{y^2z}\linedr(-1,3)
\labelL(1,3)[0]{xz^2}\linedl(1,3)\circle(1,3)
\labelL(0,2)[0]{yz^2}\lined(0,2)
\labelL(0,1)[0]{z^3}
\end{texdraw}\end{center}
\caption{$(x^2,xy,y^2)$ in degree 2, 3, and $(x,y^3)$ in degree 3}
\label{fig: 3 points in the plane}
\end{figure}

For example there are two Borel-fixed points on the 
Hilbert scheme of 3 points in the plane.  They are 
described by the Borel-fixed ideals 
$(x^2,xy,y^2)$ and $(x,y^3).$  The first one is defined in 
degree 2, the second in degree 3.  Figure 
\ref{fig: 3 points in the plane} shows these ideals, the 
first in both degrees 2 and 3, the second in degree 3. 
They are represented as filters in the posets, with the filter 
elements circled.
 
The following proposition shows $\poset(m,n)$ 
is a well-known poset.  
Let 
$\mathbf{k}$ denote the $k-$element chain on 
$\{1, 2, \ldots, k\},$ and $J(X)$ the poset on the 
order-ideals of the poset $X,$ and finally 
$\langle z_1, \ldots, z_k \rangle$ the order-ideal
generated by $z_1, \ldots, z_k.$

\begin{proposition}
We have
\[ 
   \poset(m,n) \cong J(\mathbf{m} \times \mathbf{n} ) . 
\]
In particular $\poset(m,n)$ is a distributive lattice.  
\end{proposition}

\begin{proof}
The isomorphism is given by 
\[ 
   x_0^{a_0} x_1^{a_1} \cdots x_n^{a_n} \longleftrightarrow
   \left\langle (a_0, n), (a_0+a_1, n-1), \ldots, 
      (a_0+a_1+\cdots+a_{n-1}, 1) \right\rangle            
\]
where we omit $(k,l)$ if $k = 0.$  
To see this is an isomorphism first note that any order-ideal 
$R$ 
of $\mathbf{m} \times \mathbf{n}$ is uniquely described in the 
form 
$\langle (k_1, 1), (k_2, 2), \ldots, (k_n, n) \rangle$  
by taking $k_i$ maximal such that $(k_i,i) \in R,$ or 
setting $k_i = 0$ if no such pair is in $R,$ and consider 
it as not occurring.    
Then note that $k_i \geq k_{i+1},$ for if this were not true 
then $(k_i+1, i) \leq (k_{i+1},i+1)$ which would imply 
$(k_i+1,i) \in R,$ contradicting the maximality of $k_i.$  
Hence for the $k_i$s there are unique $a_i$s such that 
$k_i = a_0+a_1+\cdots+a_{n-i-1}.$  
Finally the covering relations correspond: if 
$R = \langle (a_0, n), (a_0+a_1, n-1), \ldots, 
      (a_0+a_1+\cdots+a_{n-1}, 1) \rangle$ then
\[ 
   \begin{array}{ccc}
   \poset(m,n) & \phantom{xxxxx} & J(\mathbf{m} \times \mathbf{n}) \\
   \phantom{xxxxx} & \phantom{xxxxx} & \phantom{xxxxx} \\
   x_0^{a_0}\cdots x_i^{a_i+1}x_{i+1}^{a_{i+1}-1}\cdots x_n^{a_n} 
   & \longleftrightarrow & R \cup \{ (a_0+\cdots+a_i+1,n-i) \}    \\
   \phantom{xxxxx} & \phantom{xxxxx} & \phantom{xxxxx} \\
   \vline & & \vline \\
   \phantom{xxxxx} & \phantom{xxxxx} & \phantom{xxxxx} \\
   x_0^{a_0} \cdots x_n^{a_n} & \longleftrightarrow &  R
   \end{array}  
\]
That $\poset(m,n)$ is a distributive lattice follows 
from the \emph{fundamental theorem for finite distributive 
lattices}.  See for example \cite{stanley}, theorem 3.4.1.
\end{proof}

\begin{corollary}
\[\poset(m,n) \cong \poset(n,m)\]
\end{corollary}

Through the maps $\poset(m,n) \rightarrow 
J(\mathbf{m}\times\mathbf{n}) \rightarrow 
J(\mathbf{n}\times\mathbf{m}) \rightarrow 
\poset(n,m)$ one can 
construct the isomorphism $\poset(m,n) \cong \poset(n,m)$ 
explicitly.  Let $\x^A$ be a degree $m$ monomial 
in the variables $x_0,\ldots,x_n.$  
Write $\x^A = x_{\a_1}\cdots x_{\a_m}$ where $\a_i \leq \a_{i+1}$ 
for $i = 0, \ldots, m.$  Set
\[ 
   b_i = \left\{\begin{array}{ll}
   n-\a_m, & i=0 \\
   \a_{m-i+1}-\a_{m-i}, & i=1,\ldots, m-1 \\
   \a_1, & i=m \end{array}\right. 
\]
If 
$B = (b_0,\ldots,b_m)$ and  $\y = (y_0,\ldots, y_m)$ then the 
isomorphism $\poset(m,n) \cong \poset(n,m)$ identifies 
$\x^A$ with $\y^B.$  For example 
\[ 
   x_0^2x_1^3x_3 = x_0x_0x_1x_1x_1x_3
   \rightarrow 
   y_0^{3-3}y_1^{3-1}y_2^{1-1}y_3^{1-1}y_4^{1-0}y_5^{0-0}y_6^0
   = y_1^2y_4   .
\]

\begin{figure}[h]
\centering
\includegraphics[width=.9in]{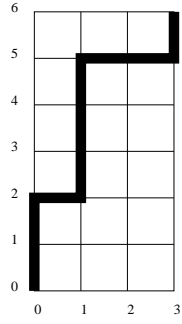}
\caption{The grid representation of the monomial 
$x_0^2x_1^3x_3 \in K[x_0,ldots,x_3].$}\label{fig: grid1}
\end{figure}

There are two much nicer ways to obtain the isomorphism.  The 
first associates to a monomial its \emph{bars and stars} 
representation.  Then one flips the role of the bars with 
that of the 
stars and reads backwards.  With the above example we find 
\[ 
   x_0^2x_1^3x_3 \sim **|***||* 
   \phantom{xxx} \longrightarrow \phantom{xxx} 
   ||*|||**| \rightarrow |**|||*|| \sim y_1^2y_4 .
\]
The second associates a monomial in $\poset(m,n)$ with a 
path in a $m \times n$ grid from the southwest corner to the 
northeast corner which always moves either up  
or to the right in integral increments.  
Each unit rise in the path indicates 
a variable corresponding to the horizontal position.  
In our example
the monomial $x_0^2x_1^3x_3$ is represented by the 
picture in figure \ref{fig: grid1}.  

\begin{figure}[h]
\centering
\includegraphics[height=.9in]{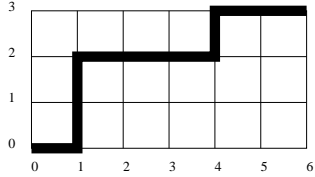}
\caption{The grid representation of the monomial 
$y_1^2y_4 \in K[y_0,\ldots,y_6]$}\label{fig: grid2}
\end{figure}

The image monomial is obtained by flipping 
the grid from southwest to northeast as in 
figure \ref{fig: grid2}.  
Once again we obtain $y_1^2y_4.$

\begin{proposition}
Let $\x^{A_1}, \x^{A_2}$ be monomials of degree 
$m$ in $K[x_0,\ldots,x_n],$ and set 
$\y^{B_1}, \y^{B_2}$ to be the corresponding monomials 
of degree $n$ in $K[y_0,\ldots,y_m]$.  Then 
\[ 
   \x^{A_1} <_{\mathrm{Lex}} \x^{A_2} \iff
   \y^{B_1} <_{\mathrm{RevLex}} \y^{B_2}   .
\]
\end{proposition}

\begin{proof}
Consider the \emph{grid representations} of the monomials as 
described above.  If $\x^{A_1} <_{\mathrm{Lex}} \x^{A_2}$
then the first juncture at which the path corresponding to $\x^{A_1}$ 
differs from that of $\x^{A_2}$ must have the latter 
path going up, while the former going right.  Hence after the 
flip the last juncture in which the paths for $\y^{B_1}$ and 
$\y^{B_2}$ differ will have the latter path coming in from 
the left, while the former from underneath (see 
figure \ref{fig: grid3}).  One readily verifies 
this is equivalent to $\y^{B_1} <_{\mathrm{RevLex}} \y^{B_2}.$  

\begin{figure}[h]
\centering
\includegraphics[width=2.5in]{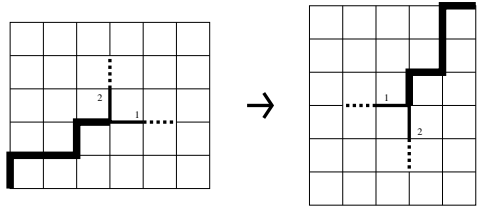}
\caption{The effect of a flip on the ordering of 
two monomials.}\label{fig: grid3}
\end{figure}

Alternatively one can use the \emph{bars and stars} 
representations and use the same logic. 
\end{proof}

We include here a lemma which we find demonstrates 
the interplay of the combinatorics of 
filters in $\poset(m,n),$ and the algebra of 
their defining ideals.  

\begin{lemma} \label{lemma: when is a monomial in I_m+1}
Let $I$ be a Borel-fixed ideal defined in 
degrees $\leq m.$  Let $\F = I \cap \poset(m,n)$ 
be the corresponding filter.  If $\x^A$ is a standard 
monomial of degree $m$ ($\x^A \in S_m \setminus I_m$), then
\[ 
   x_i \x^A \in I \iff \frac{x_i}{x_{\max(A)}}\x^A \in \F. 
\]
In particular, if $\x^A$ is \emph{Borel maximal} in 
$\poset(m,n) \setminus \F$ (that is every greater monomial 
lies in $\F$) then 
\[
	x_i \x^A \in I \iff i < \max(A).
\]
\end{lemma}

\begin{proof}
Let $k = \max(A).$  If $x_i/x_k \x^A \in \F$ then 
certainly $x_i \x^A \in I.$  Conversely suppose 
$x_i\x^A \in I.$  Since $I$ is generated in degrees less 
than or equal to $m,$ there must be a monomial 
$\x^B \in \F$ and a variable $x_j$ such that 
\[ 
   x_i \x^A = x_j \x^B .
\]
Certainly $i \neq j$ (since $\x^A \notin I$) and hence 
$j \leq \max(A) = k.$  But then
\[ 
   \frac{x_i}{x_k}\x^A = \frac{x_j}{x_k}\x^B \in \F 
\]
since $\F$ is a filter.
\end{proof}


\section{The tangent space to the Hilbert scheme at a Borel-fixed point.}
\label{Description of the tangent space}

Fix a projective space $\P^n$ and a Hilbert polynomial 
$\poly(z)$ and set $\Hilb$ to be the corresponding Hilbert scheme.  
If $z \in \Hilb$ corresponds to the scheme $Z \in \P^n$ then it is 
known that 
\[ 
   T_z \Hilb = \HH^0 \N_{Z/\P^n} .
\]
That is, the tangent space to $z \in \Hilb$ is identified with 
the global sections of the normal sheaf to $Z \subset \P^n.$  

Let $z$ be a Borel-fixed point on $\Hilb,$ let the corresponding ideal 
sheaf be $\I,$ and let $m$ be the Gotzmann number for the corresponding 
Hilbert polynomial.  The tangent space to this point 
will be a subspace of the tangent space to the Grassmanian 
$\Grass(r,S_m).$  Recall what the tangent space to the Grassmanian 
looks like.  Let $v \in \Grass(r,S_m)$ correspond to 
the $r-$dimensional subspace $V$ of $S_m.$  Then
\[ 
   T_v \Grass(r,S_m) = \Hom \left( V, S_m / V ) \right). 
\]
We want to see how the tangent space to $z \in \Hilb$ sits naturally
as a subspace of $T_{z} \Grass(r,P_m)$. 
Recall that $\N_{Z/\P^n} = \sheafHom (\I / \I^2, \O_Z)$
so 
$T_z \Hilb = \HH^0 \N_{Z/\P^n} = 
\Hom (\I / \I^2, \O_Z)$.  Now a map
$\phi : \I / \I^2 \rightarrow \O_Z$ composes with the natural
map $\I \rightarrow \I / \I^2$ to give 
$\tilde{\phi} : \I \rightarrow \O_Z.$  Twisting by $m$ we then
get a map $ \I (m) \rightarrow \O_Z(m).$  Then take global sections
to get $ \HH^0 ( \I  (m) ) \rightarrow \HH^0 ( \O_Z(m) ) $.  
Let $I = \bigoplus_{d \geq0} \HH^0(\I(d)).$  Since 
$m \geq \reg(\I),$ we see 
\[ 
   \HH^0 ( \I (m) ) = I_m \mathrm{\ and\ }
   \HH^0 (\O_Z(m)) = S_m/I_m .
\]  
Hence we get the map
$I_m \rightarrow S_m/I_m \in T_{z} \Grass(r,S_m).$ 

The tangent space to a point $z \in \Hilb$ can also be identified 
with the space of first-order infinitesimal deformations of 
the corresponding scheme in $\P^n.$  Specifically, if $z \in \Hilb$ 
corresponds to the scheme defined by the ideal $I = I^{\sat}$ then 
\begin{eqnarray*}
   T_z \Hilb & \cong & \{ J \subset S[\e] \mid J \mathrm{\ flat\ over\ } 
                K[\e],\ J|_{\e = 0} = I \} \\
          & \cong & \{ J_m \subset S_m[\e] \mid J_{\geq m} 
          \mathrm{\ flat\ over\ } 
                K[\e],\ J_{\geq m}|_{\e = 0} = I_{\geq m} \}
\end{eqnarray*}

Now let $z \in \Hilb$ be a Borel-fixed point whose corresponding 
scheme is defined by the Borel-fixed ideal $I = I^{\sat}.$  
Set $\F= I\cap \poset(m,n)$ to be the corresponding 
filter, and $\R = \poset(m,n)\setminus \F$ 
the corresponding order-ideal of the standard 
monomials in degree $m,$ as we defined in section 
\ref{The poset P(m,n)}.  
For notational reasons it is often 
simpler if we consider $\F$ and $\R$ as consisting of the 
exponent vectors of the monomials.  In this 
respect we will switch 
back and forth between monomials and their exponent 
vectors and trust that no confusion will arise.  
Note that since the Borel subgroup fixes $z$ it 
induces an action on $T_z \Hilb.$  

An arbitrary vector in the tangent space to 
$z$ in $\Grass(r,S_m)$ is 
given by a $K-$linear map 
$ \phi : I_m \rightarrow S_m/I_m $ which can be 
described uniquely by 
\[ 
   \phi(\x^A) = \sum_{B \in \R} c_{AB} \x^B, \quad A \in \F.
\]
We denote this by the doubly-indexed vector  
\[ 
   \left( c_{AB} \right)_{A\in\F, B\in\R} 
   = \left( c_{AB} \right) \in K^{\#\F \times \R} 
   \cong T_z \Grass(r,S_m). 
\]
This vector then lies in $T_z \Hilb$ if and only if 
\[ 
   J := \left( \x^A + \e \sum_{B \in \R} c_{AB} \x^B \mid A \in \F 
        \right) ,
\]
with $\e^2=0,$ is flat over $K[\e]$ (that 
is $J$ defines a first order infintesimal deformation of $I$).  
We will set 
$\{e_{AB} \mid A \in \F, B \in \R \}$ to be the basis  
of $T_z\Grass(r, S_m)$ where $e_{AB}$ is the vector 
with $0$s in every coordinate except the one with index $(A,B).$  


\section{The maximal torus eigenvectors}
\label{The maximal torus eigenvectors}

In this section we classify those vectors of the 
tangent space to the Hilbert scheme at Borel-fixed 
point which are eigenvectors for the maximal torus 
subgroup of $\GL(n+1,K)$ consisting of diagonal matrices.  

\begin{lemma} \label{the torus-fixed directions}
Let $z$ be a Borel-fixed point on $\Hilb$.  Let $m$ be the regularity of 
the sheaf of ideals $\I$ of $Z$.  Let $I$ be the 
(Borel-fixed) ideal given by 
\[ 
   I_{d} = \left\{
   \begin{array}{lr}
   H^0 \I(d) & \mathrm{\ if\ } d \geq m \\
   0          & \mathrm{\ otherwise}
   \end{array} \right.
\]
Set $\F$ to be the set of (exponent vectors of) monomials of $I$ of 
degree $m.$  Set $\R$ to be the set of (exponent vectors of) 
the standard monomials of degree $m.$ 
(Recall that $\F$ is a filter of $\poset(m,n),$ while 
$\R$ is the complimentary order-ideal.)  Then the infinitesimal deformation 
\[ 
   \left(\x^A + \e\sum_{B\in \R} c_{AB} \x^B \mid A \in \F\right) 
   = (c_{AB})_{A\in \F, B\in \R}  
\]
considered as an element of the tangent space to $z$ is an 
eigenvector for the maximal torus of diagonal matrices if and only if
there exists $K \in \Z^{n+1}$ of degree $0$ such that 
\[ 
   c_{AB} \not = 0 \implies B-A = K .
\]
\end{lemma}

\begin{proof}
Let $\Lambda$ be the diagonal matrix with diagonal entries 
$\lambda = (\lambda_0, \ldots, \lambda_n)$.  For the deformation
$(c_{AB})_{A\in \F, B\in\R}$ set 
\[ 
   r_A = \x^A + \e\sum_{B\in \R} c_{AB} \x^B .
\]
Then 
\[ 
   \Lambda \cdot r_A = \lambda^A\x^A + \e\sum_{B\in \R} c_{AB} \lambda^B\x^B 
\]
from which we get 
\[ 
   \lambda^{-A} \Lambda \cdot r_A = 
   \x^A + \e\sum_{B\in \R} c_{AB} \lambda^{B-A}\x^B 
\]
Hence 
\[ 
   \Lambda \cdot (c_{AB}) = (\lambda^{B-A} c_{AB}) 
\]
By varying the $\lambda_i$'s we see that the $B-A$ must be constant 
over those $A,B$ for which $c_{AB} \not = 0.$  Setting 
$K = B-A$ we get our result, and in fact 
\[  
   \Lambda \cdot (c_{AB}) = \lambda^{K} (c_{AB}) 
\] 
\end{proof}


\section{The Borel eigenvectors}\label{The Borel eigenvectors}

For an eigenvector of the maximal torus as in lemma 
\ref{the torus-fixed directions} set  
\[
   \F' = \F \setminus \{ A \in \F \mid \exists B \in \R \mathrm{\ such\ that\ } 
          c_{AB} \neq 0 \} 
\]
and
\[ 
   \F'' = \F \cup \{B \in \R \mid \exists A \in \F \mathrm{\ such\ that\ } 
          c_{AB} \neq 0 \} 
\]
We will say this eigenvector has \emph{type} $(\F', \F'', K).$  
Note that 
$\F''$ is determined in terms of $\F'$ and $K$: 
\[
   \F'' = \F \cup \left( (\F \setminus \F') + K \right) 
   \quad\mathrm{(disjoint\ union)}.
\]
So we may refer to this vector as having type $(\F', K).$  

If $(c_{AB})_{A\in\F,B\in\R}$ is an eigenvector for the maximal torus 
we can denote it by $(c_A)_{A\in\F}$ without 
confusion since for any $A\in\F$ there 
is at most one $B \in \R$ (namely $A+K$) 
such that $c_{AB} \neq 0.$  Specifically, the notation 
$(c_A)_{A \in \F}$ (or even more simply $(c_A)$) for 
an eigenvector of type $(\F', K)$ will refer 
to the ideal of $S[\e]$ generated by the elements 
\[
	\x^A, \mathrm{\ for\ each\ } A \in \F \setminus \F', \quad 
		\x^A + \e c_A \x^{A+K}, \mathrm{\ for\ each\ } A \in \F'. 
\]

The question of which vectors in $T_z \Hilb$ are eigenvectors 
for the Borel group of upper-triangular matrices is a little 
more tricky.  First observe that since the diagonal matrices 
are a subgroup of the Borel group we must have that such a 
vector is an eigenvector for the maximal torus.  

Take an eigenvector for the maximal torus, of 
type $(\F', \F'', K).$  We say this vector is a 
{\it pseudo-eigenvector} 
for the Borel subgroup if its image under the action 
of any upper-triangular matrix also has type $(\F',\F'', K).$  
Certainly an eigenvector for the Borel group is a
pseudo-eigenvector. 

Let $E_i \in \Z^{n+1}, i = 0,\ldots,n$ be the vector with a $1$ in the 
$i$'th position, and $0$'s elsewhere (note that the `0 position' is the 
first coordinate).  For $i = 1,\ldots,n$ set 
$\D_i = E_{i-1}-E_i.$  Notice that every covering relation in 
$\poset(m,n)$ is 
of the form 
$A \prec_{\mathrm{B}} A+\D_i$ for some $i.$  

\begin{lemma} \label{almost Borel-fixed = Borel-fixed}
Let $(c_A)_{A\in\F}$ be an eigenvector for the maximal torus, of 
type $(\F', \F'', K).$  If $(c_A)$ is a pseudo-eigenvector 
for the Borel subgroup then:
\begin{itemize}
\item[(i)]  $\F'$ and $\F''$ are both filters;
\item[(ii)]  If $A, A+\Delta_i \in \F \setminus \F'$ then 
\[ 
   c_{A+\Delta_i} = \frac{a_i}{b_i} c_A 
\]
where $a_i = \deg_i(A)$ and $b_i = deg_i(A+K).$  
\item[(iii)]  $(c_A)_{A\in\F}$ is an eigenvector for the Borel 
subgroup.  
\end{itemize}
\end{lemma}

\begin{proof}
Let $J$ be the corresponding ideal of $S[\e].$  Its generators 
are 
\[ 
   \x^A, A \in \F', \quad \x^A + \e c_A \x^{A+K}, A \in \F \setminus \F' .
\]
Let $h$ be an arbitrary 
upper-triangular matrix.  By assumption $hJ$ has generators 
$\x^A + \e s_A \x^{A+K}$ where $s_A = 0 \iff A \in \F'.$  
Hence for $B \in \F$ there must be a relation
\[ 
   \x^{B} + \e s_{B} \x^{B+K} = \sum_{A \in \F'} (\lambda_A+\e\mu_A) 
   h \cdot \x^A + \sum_{A \in \F \setminus \F' }
   (\lambda_A+\e\mu_A) h \cdot 
   (\x^A + \e c_A \x^{A+K}) .
\]
This is equivalent to the two equations
\begin{equation} \label{eqA'}
\x^{B} = h \left( \sum_{A \in \F} \lambda_A \x^A \right)
\end{equation}
and
\begin{equation} \label{eqB'}
   s_{B} \x^{B+K} = h \left(
   \sum_{A \in \F} \mu_A\x^A + 
   \sum_{A \in \F \setminus \F'} \lambda_A c_A x^{A+K} \right)
\end{equation}
If $B \in \F'$ then $s_{B} = 0,$ and (\ref{eqB'}) becomes 
\[ 
   0 = h \left( 
   \sum_{A \in \F} \mu_A \x^A + 
   \sum_{A \in \F \setminus \F'} \lambda_A c_A \x^{A+K} 
   \right)   
\]
Since $h$ is nonsingular and the sum on the right is over 
distinct monomials, we find $\mu_A = 0, A \in \F$ and 
$\lambda_A = 0, A \in \F \setminus \F'.$  So (\ref{eqA'}) becomes 
\[ 
   h^{-1} \x^{B} = \sum_{A \in \F'} \lambda_A \x^A .
\]
Since $B \in \F'$ is arbitrary, as well as $h,$ we see 
$\F'$ is a filter.

On the other hand if $B \in \F \setminus \F'$ then $s_{B} \neq 0$ and 
(\ref{eqA'}) and (\ref{eqB'}) can be rewritten 
\[ 
   h^{-1}\x^{B} = \sum_{A \in \F} \lambda_A \x^A 
\] 
and
\[ 
   h^{-1} \x^{B+K} = 
   \sum_{A \in \F} s_{B}^{-1} \mu_A \x^A + 
   \sum_{A \in \F \setminus \F'} s_{B}^{-1} \lambda_A c_A \x^{A+K} 
\]
Since $\F'' = \F \du ( (\F \setminus \F') + K )$ 
these show $\F''$ is a filter.  

Lastly let $A \in \F \setminus \F'.$  
Set $a_i = \deg_i(A),$ and $b_i = \deg_i(A+K).$ 

Let $h_i$ be the upper-triangular 
matrix which sends $x_i \mapsto x_i+x_{i-1}$ and 
leaves the other variables fixed.  We compute
\[ 
   h_i \left( \x^A + \e c_A \x^{A+K} \right) = 
   \sum_{j=0}^{a_i} {a_i\choose j} \x^{A+j\D_i} + 
   \sum_{j=0}^{b_i} {b_i\choose j} c_A \e x^{A+K+j\D_i}   
\]
where $a_i = \deg_i(A),$ and $b_i = \deg_i(A+K).$  
Set $l$ to the maximal such that $A+K+l\D_i \in \R.$  We already know that 
$\F''$ is a filter, and this implies that $l \leq a_i.$  
Then $h_iJ$ contains 
\[ 
   h_i \left( \x^A + \e c_A \x^{A+K} \right)  
   = \sum_{j=l+1}^{a_i} {a_i \choose j} 
   \underbrace{\x^{A+j\D_i}}_{\in h_iJ} + 
   \sum_{j=0}^{l} {a_i \choose j} 
   \left[ \x^{A+j\D_i} + \e c_A \frac{{b_i\choose j}}{{a_i\choose j}} 
   \x^{A+K+j\D_i} \right]   
\]
Each individual term in the left-hand sum lies in $h_iJ$ so the right-hand 
sum must be in $h_iJ.$  By assumption $h_i J$ has generators of the 
form 
\[ 
   \{ \x^B + \e s_B \x^{B+K} \mid B \in \F \}. 
\]
So
$\x^{A+j\D_i}+\e s_{A+j\D_i}\x^{A+K+j\D_i} \in h_iJ.$  It follows that 
\[ 
   \e \sum_{j=0}^{l} {a_i \choose j} \left( 
   s_{A+j\D_i} -  c_A \frac{{b_i\choose j}}{{a_i\choose j}} 
   \right) \x^{A+K+j\D_i} \in h_iJ . 
\]
But since $h_i J$ is flat over $K[\e]$ (recall that the Borel group 
acts on $T_z\Hilb$) this can only happen 
if all the coefficients are 
zero.  For $j=0$ this says $s_A = c_A,$ and this holds for any 
$A \in \F \setminus \F'.$  This gives (iii) since the $h_i$ 
together with the diagonal matrices generate the Borel subgroup.  
If $A+\D_i \in \F \setminus \F'$ as 
well then $l \geq 1.$  Setting the $j=1$ coefficient to be zero 
gives us $c_{A+\Delta_i} = \frac{b_i}{a_i} c_A.$
\end{proof}

This gives us the following description of the Borel-eigenvectors of 
the tangent space:

\begin{theorem} \label{the borel-fixed directions} 
Let $z$ be a Borel-fixed point on $\Hilb.$  
Let $m$ be the regularity of the
sheaf of ideals $\I$ defining the subscheme corresponding to $z.$ 
Let $I$ be the ideal given by 
\[ 
   I_{d} = \left\{
   \begin{array}{lr}
   H^0 \I(d) & \mathrm{\ if\ } d \geq m \\
   0          & \mathrm{\ otherwise}
   \end{array} \right.
\]
Set $\F$ to be the set of (exponent vectors of) monomials of 
$I$ of degree $m$ (which is a filter of $\poset(m,n)$) and 
$\R$ the set of (exponent vectors of) the degree $m$ standard monomials.  
Then the infinitesimal deformation 
\[ 
   \left(\x^A + \e\sum_{B\in \R} c_{AB} \x^B \mid A \in \F\right) 
   = (c_{AB})_{A\in F, B\in \R}  
\]
is Borel-fixed as an element of the tangent space to $Z$ 
if and only if
\begin{enumerate}
\item There exists $K \in \Z^{n+1}$ of degree $0$ such that 
      $ c_{AB} \not = 0 \implies B-A = K .$
\item let 
      \[ 
         \F' = \F \setminus \{ A \in \F \mid \exists B \in \R, c_{AB} \neq 0\} 
      \]
      and 
      \[ 
         \F'' = \F \cup \{ B \in \R \mid \exists A \in \F, c_{AB} \neq 0 \} 
\]
Then 
\begin{itemize}
\item[(i)] $\F'$ and $\F''$ are filters; 
\item[(ii)] For each $A = (a_0, \ldots, a_n) \in \F \setminus \F'$, 
$B = A + K = (b_0, \ldots, b_n)$, if $A + \D_i \in \F \setminus \F'$ then 
      $c_{A+\D_i,B+\D_i} = \frac{b_i}{a_i}c_{AB}.$
\end{itemize}
\end{enumerate}
\end{theorem}

\begin{figure}[h]
\begin{center}
\includegraphics[width=2in]{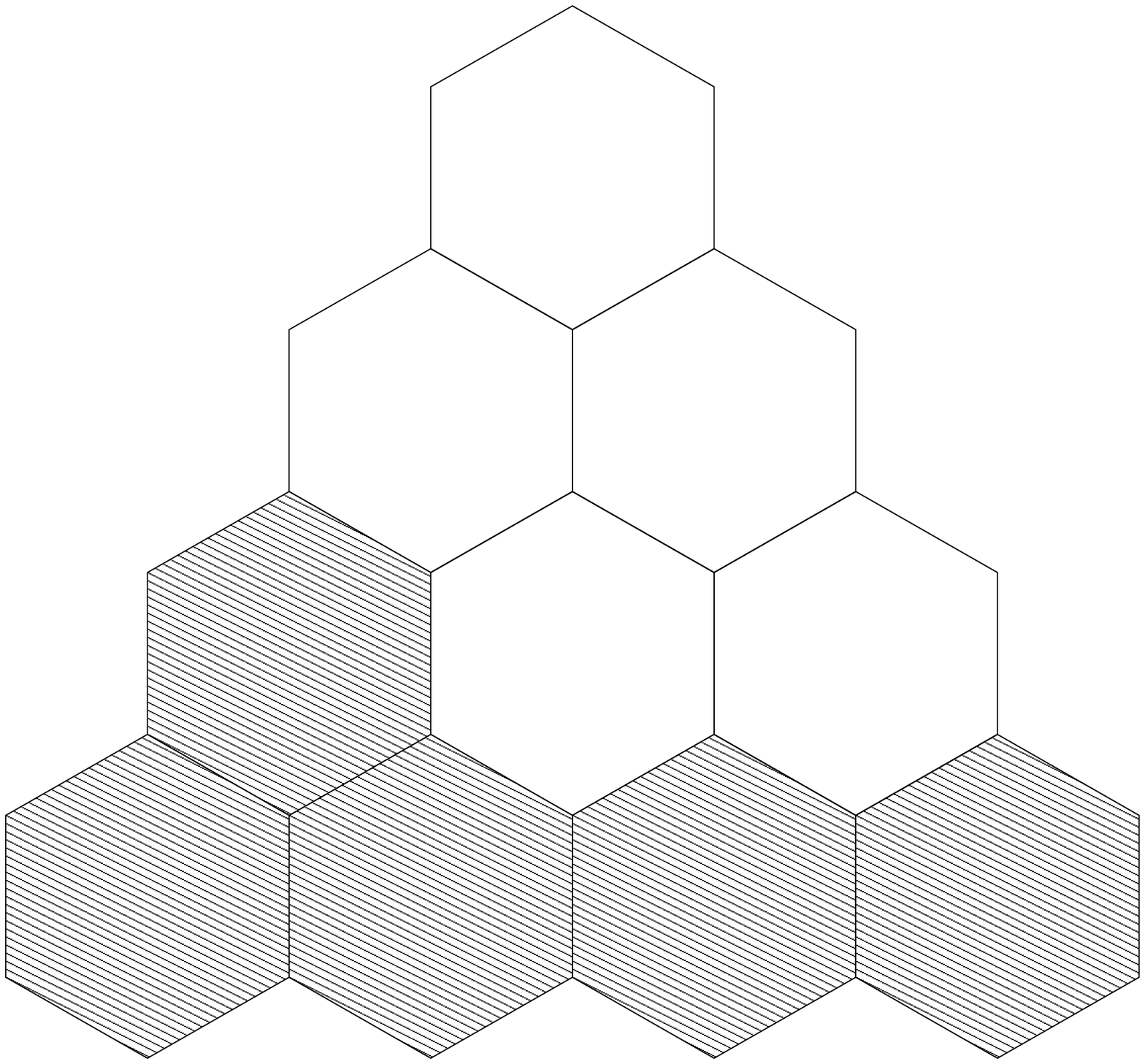}
\put(-130, 15){$x^3$}
\put(-95, 15){$x^2y$}
\put(-60, 15){$xy^2$}
\put(-25, 15){$y^3$}
\put(-115, 45){$x^2z$}\put(-80, 45){$xyz$}\put(-45, 45){$y^2z$}
\put(-95, 80){$xz^2$}\put(-60,80){$yz^2$}\put(-75, 110){$z^3$}
\caption{The monomials of degree 3.}
\label{fig: bf0}
\end{center}
\end{figure}

The characterization of the Borel-fixed tangent vectors in 
\ref{the borel-fixed directions} is rather dry and 
unenlightening, so let us illustrate with some examples.  
Take the Borel-fixed ideal 
$I = (x^3, x^2y, xy^2, y^3, x^2z) \subseteq K[x,y,z].$  
We depict the monomials in degree $3$ as in figure 
\ref{fig: bf0}, with the monomials in $I$ shaded.

\begin{figure}[h]
\begin{center}
\includegraphics[width=2in]{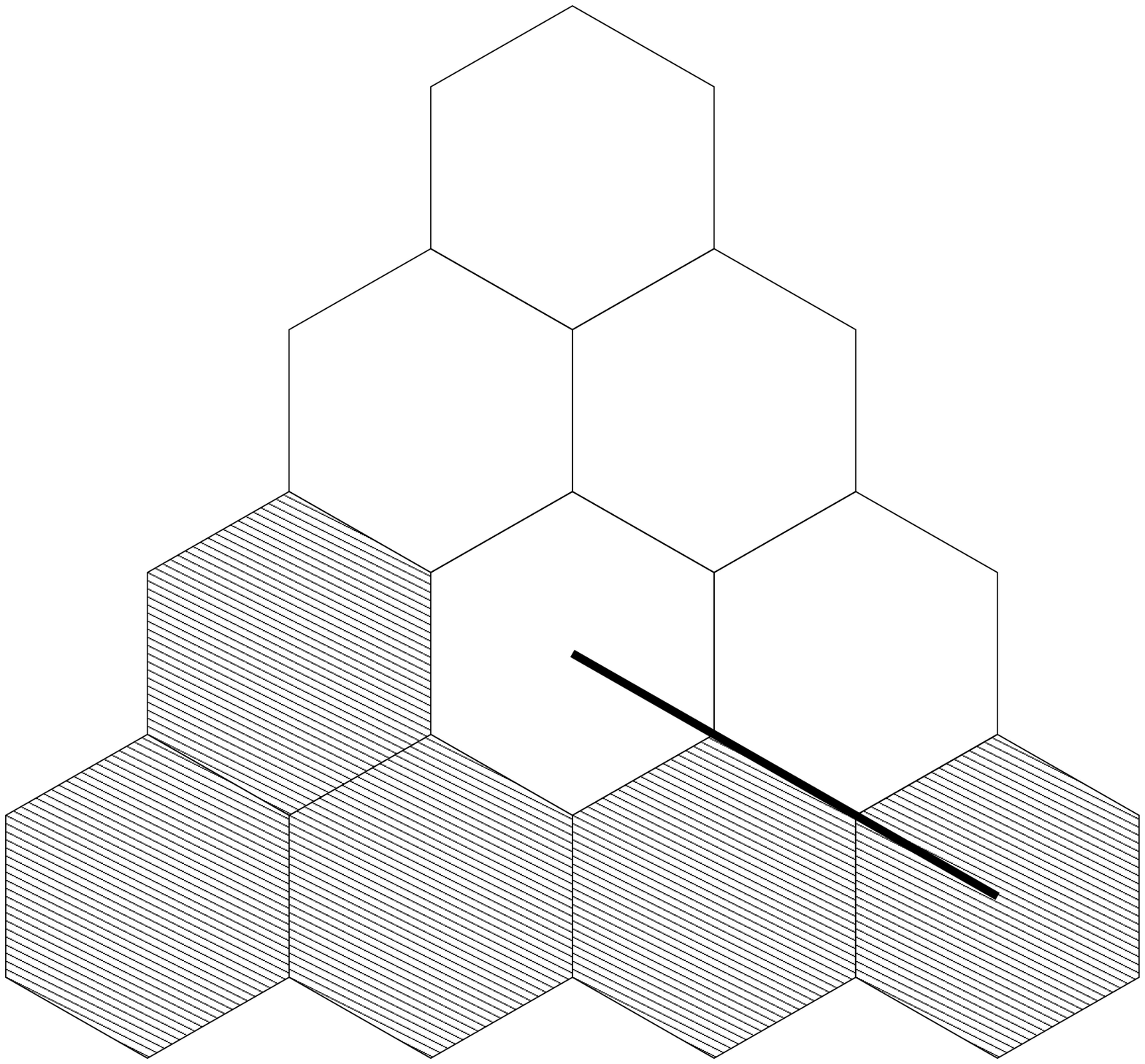}
\quad
\includegraphics[width=2in]{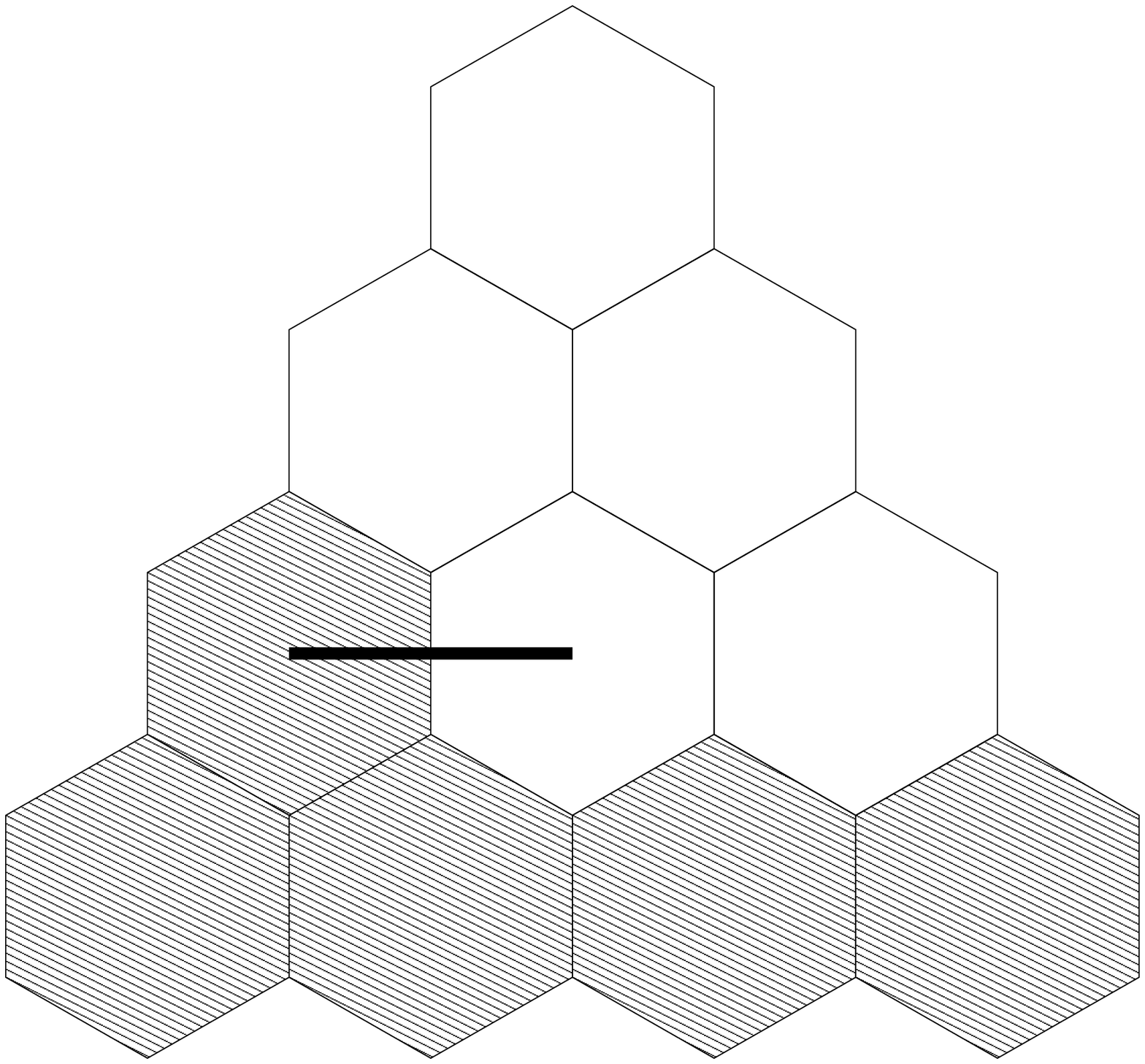}
\caption{Two Borel eigenvectors.}
\label{fig: bf1}
\end{center}
\end{figure}

We will represent a 
tangent vector $(c_{AB})$ by placing line segments 
on this picture, with a line segment extending from 
the hexagon representing $A$ to that representing $B$ if 
$c_{AB} \neq 0,$ and labeling that line with the number 
$c_{AB}$ if it is not $1.$  So for instance, on the left
of figure 
\ref{fig: bf1} we depict the vector corresponding to the 
infintesimally deformed ideal 
$(x^3, x^2y, xy^2, y^3 + \e xyz, x^2z) \subseteq K[x,y,z][\e],$ 
where $\e^2 = 0.$  Condition 1 of theorem 
\ref{the borel-fixed directions} is automatically 
satisfied in this example.  Furthermore, here we 
have 
$\F' = \{ x^3, x^2y, xy^2, x^2z \},$ and 
$\F'' = \{ x^3, x^2y, xy^2, y^3, x^2z, xyz \},$ which 
are both filters of $\poset(3,2),$ and hence condition 
2(i) is satisfied.  Finally 2(ii) is immediate.  Hence 
this vector is an eigenvector for the Borel subgroup.  
One can easily check the vector depicted on the 
right of figure \ref{fig: bf1} also is a Borel eigenvector.  

\begin{figure}[h]
\begin{center}
\includegraphics[width=2in]{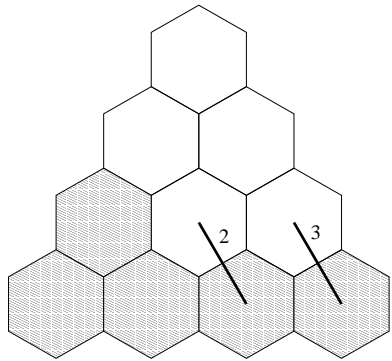}
\caption{Another Borel eigenvector.}
\label{fig: bf2}
\end{center}
\end{figure}

Now consider instead the tangent vector depicted in 
figure \ref{fig: bf2}.  This represents the 
ideal $(x^3, x^2y, xy^2+2\e xyz, y^3 + 3\e y^2z, x^2z) 
\subseteq K[x,y,z][\e].$  Condition 1 of the theorem 
is satisfied with $K = (-1,1,0).$  In fact one can 
see that condition 1 just requires that all the line 
segments that appear must be rigid translates of 
each other.  Condition 2(i) can be easily verified 
(note that a set of monomials in this picture form a 
filter if and only if they are closed under taking 
steps down and steps left).  Finally condition 2(ii) holds:  
in the notation of the theorem we have $a_i = 3$ 
(the degree of the $y$ variable in $y^3$) and 
$b_i = 2$ (the degree of the $y$ variable in $y^2z$).  
So this vector is also a Borel eigenvector.  

By examining possible Borel eigenvectors to be placed 
on figure \ref{fig: bf0} one can see that the three 
depicted in figures \ref{fig: bf1} and \ref{fig: bf2} 
are the only eigenvectors (up to scalar multiples of 
course).  

We conclude this chapter by giving an alternative way of 
viewing theorem \ref{the borel-fixed directions}.  Note that 
it gives a characterization of those lines through the 
origin of the tangent space of a Borel-fixed point which are 
fixed by the action of the Borel subgroup of $\GL(n+1,K).$  
However the lines through the origin in the tangent space 
to any point on any scheme can be viewed as points in 
the blowup of the scheme at that point.  These are 
called \emph{infintely near points}.  In this language 
theorem 
\ref{the borel-fixed directions} characterizes 
\emph{the infintely near 
Borel-fixed points on the Hilbert scheme}.  


\section{An infinitesimal version of Galligo's theorem}
\label{An infinitesimal version of Galligo's theorem}

In this section we prove the main result of this paper, namely 
that an infinitesimal version of Galligo's theorem 
\cite{galligo} holds.  First let us recall Galligo's result.
Recall that $\GL_{n+1}(K)$ acts on the set of ideals 
of the polynomial ring over $K$ with $n+1$ variables. 

\begin{theorem}[Galligo]  Let 
$I \subseteq K[x_0, x_1, \ldots, x_n]$ be any ideal and 
let $>$ be a monomial order with, say, 
$x_0 > x_1 > \ldots > x_n.$    
Then there exists a Zariski open (and therefore dense) 
subset $U$ of $\GL_{n+1}(K)$ such that for $g \in U$ the 
initial ideal 
$\init_> (gI)$ is constant over $g \in U$ and Borel-fixed 
(that is fixed by the action of the subgroup of 
upper-triangular matrices).  
\end{theorem}

In order to state our infinitesimal version we first need to 
switch from monomials orders $>$ to weight vectors.  
A {\it weight vector} $w$ is an element of $\Z^{n+1}$ which 
we use to partially order monomials by associating associating 
to each monomial a {\it weight} $w(\x^{A}) = w \cdot A$:    
\[
	\x^A >_w \x^B \iff w\cdot A > w\cdot B.
\]
Note that it makes perfect sense to consider weight vectors 
as any element of $\mathbf{R}^{n+1}.$  However 
we will restrict our weight vectors 
to have integer coordinates.  

If $>$ is a monomial order, and $m$ is a positive integer, then 
we will say that a weight vector $w$ {\it induces $>$ in degree $m$} 
if for any monomials $\x^A, \x^B$ of degree $m$ we have 
$\x^A > \x^B \iff \x^A >_w \x^B.$  If $w$ induces $>$ in degree $m$ then 
$w$ induces $>$ in every degree $m' \leq m.$   
For any monomial order $>$ and degree $m$ there is a weight 
vector (in fact many!) which induce $>$ in degree $m.$  
The weight vector $w$ will be said 
to {\it distinguish monomials in degree $m$} if for 
any two distinct monomials $\x^A, \x^B$ of degree $m$ we have 
either $\x^A >_w \x^B$ or $\x^A <_w \x^B.$  If $w$ distinguishes 
monomials in degree $m$ then $w$ distinguishes monomials in 
every degree $m' \leq m.$  

Set $S = K[x_0, \ldots, x_n]$ and let $I \subseteq S$ be an ideal.  
If $w$ is a weight vector then for $t \neq 0$ we set 
$\l = \l(t)$ to be the diagonal matrix with diagonal entries 
$t^{-w_0}, \ldots, t^{-w_n}.$  We then define the new 
ideal $I(t)$ by $I(t) = \l(t) \cdot I.$  This is the result of 
replacing every occurrence of the variable $x_i$ in $I$ by 
$t^{-w_i}x_i,$ for each $i = 0, \ldots, n.$  Then the 
set $\{ I(t) \mid t \neq 0\}$ forms a one-parameter family 
of ideals which we can view as a curve on the Hilbert scheme.  
We call limit ideal the {\it initial ideal with respect to $w$}: 
\[
	\init_w (I) := \lim_{t \rightarrow 0} I(t).
\]
If $w$ induces the monomial order $>$ to a large enough degree (the 
Gotzmann number is plenty large enough) then 
$\init_w (I) = \init_> (I).$  

The goal of this section is to prove the following infinitesimal version 
of Galligo's theorem:  If $I \subseteq S$ is an ideal, and if 
$w$ is a weight vector, with say $w_0 > w_1 > \ldots > w_n,$ 
which distinguishes monomials in a large 
enough degree, then for $g$ in a 
dense open subset of $\GL_{n+1}(K)$ the family of ideals 
$\{ (gI)(t) \}$ (as defined above) has for limit 
$\init_w (gI)$ a Borel-fixed ideal (Galligo's theorem), and the 
direction this family, viewed as lying on the Hilbert scheme, 
approaches this limit is itself Borel-fixed.  We remark that 
since the Borel group fixes the limit point on the Hilbert scheme, 
it descends to an action on its tangent space.  

Fix a projective space $\P^n$ with homogeneous coordinate 
ring $S = K [ x_0, \ldots, x_n ].$  For $r > 0$ and 
$m > 0$ consider the vector space $\wedge^r S_m$ of 
$r-$fold wedge products of degree $m$ homogeneous forms in $S.$  
In this space we call the wedge product of $r$ pair-wise 
distinct monomials 
$\x^{B_1}\wedge\cdots\wedge\x^{B_r}$ a \emph{state}. 
Two states 
$\wedge\x^{B_i} = \x^{B_1}\wedge\cdots\wedge\x^{B_r}$ 
and 
$\wedge\x^{C_i} = \x^{C_1}\wedge\cdots\wedge\x^{C_r}$ 
span the same linear (one dimensional) 
subspace iff there is a 
permutation $\sigma \in \mathrm{Sym}(r)$ such that 
$B_i = C_{\sigma (i)},$ in which case 
$\wedge\x^{B_i} = (-1)^{\sigma}(\wedge\x^{C_i}).$  
Two such states will be called \emph{equivalent}.   
The {\it associated monomial} of a state 
$\x^{B_1}\wedge\cdots\wedge\x^{B_r}$ 
is the degree $rm$ monomial 
$\x^{B_1}\cdots\x^{B_r} = \x^{B_1 + \cdots + B_r}.$  
If we are given a weight vector then we declare 
the \emph{weight} of a state as the weight of its 
associated monomial.  
Equivalent states have the same weight.  However 
nonequivalent states may still have the same weight;  
for instance $xz \wedge y^2$ and $xy \wedge yz$ each 
have weight will have the same weight for {\it any} 
weight vector as they both have the same associated 
monomial $xy^2z.$

Every element $f \in \wedge^r S_m$ can be written 
uniquely as a linear combination of states, if one ignores 
the distinction of equivalent weights.  We define 
the \emph{support} of $f,$ denoted $\supp(f)$ to be 
the set of those states appearing with non-zero 
coefficients.  If we have a weight vector then 
for a given weight value $N$ (which is 
an integer) we define 
$\supp_{N}(f)$ to be the set of those states of 
weight $N.$  

The individual summands in the expression of 
$f \in \wedge^rS_m$
as a linear combination of states will be referred to 
as the \emph{terms} of $f.$   Given a weight vector, the 
weight of such a term 
is just the weight of the associated monomial. 
For a given weight $N$ we will write 
$f_N$ for the sum of the terms of $f$ with weight $N.$  
We will write 
$f_{\geq N}$ for $\sum_{M \geq N} f_A,$ the sum of the 
terms of $f$ with weight at least $N.$  We 
analogously define $f_{>N}, f_{<N},$ and $f_{\leq N}.$  
We will also set 
$\supp_{\geq N}(f) = \supp(f_{\geq N}),$ and 
similarly for $>, <, \leq.$  

Given  linearly independent homogeneous forms 
$f_1, \ldots, f_r$ of degree $m,$ and a weight vector $w$ 
that distinguishes monomials of degree $m,$ we note that 
$f = f_1\wedge\cdots\wedge f_r$ has a unique 
term of maximal weight;  namely pick 
$f_1', \ldots, f_r'$ to span the same subspace as 
$f_1, \ldots, f_r,$ and such that the initial term 
(that is the term of largest weight) of each 
$f_i'$ does not appear in any other $f_j'.$  Then 
$f_1'\wedge\cdots\wedge f_r'$ differs from 
$f_1\wedge\cdots\wedge f_r$ only by a non-zero scalar 
multiple and it has a unique highest term 
of the form 
$c \cdot \init_w(f_1')\wedge\cdots\wedge \init_w(f_r'),$ 
with $c\neq 0.$   
We will write $\init_w (f)$ for the corresponding 
state 
$\init_w(f_1')\wedge\cdots\wedge \init_w(f_r') \in \supp (f)$ 
of highest weight.  

Let $G$ be an $(n+1) \times (n+1)$ matrix with 
variable entries $G_{ij}.$  If 
$f = f_1 \wedge \cdots \wedge f_r \in \wedge^r S_m,$
we can express 
$Gf = (Gf_1) \wedge \cdots (Gf_r)$ as a linear combination 
of states whose coefficients are polynomials in the 
variables $G_{ij}.$  Let $U(f) = U \subseteq \GL(n+1,K)$ be the 
open set where none of the non-identically zero polynomials 
vanish.  We make the following observations which follow immediately:
\begin{itemize}
\item[(1)]  If $g, g' \in U = U(f)$ 
	then $\supp(gf) = \supp(g'f);$
\item[(2)]  If $g \in U$ and $g' \in \GL_{n+1}$ is 
	arbitrary then $\supp(g'f) \subseteq \supp(gf).$  
\item[(3)]  If $\lambda \in \GL_{n+1}$ is a diagonal matrix than 
	$\lambda U = U$
\end{itemize}

In addition we get the following:

\begin{lemma} \label{support equality}
Let $f = f_1 \wedge \cdots \wedge f_r \in \wedge^r S_m,$ and 
let $w$ be a weight vector which distinguishes monomials in 
degree $rm.$  Set $\l = \l_w(t)$ be the diagonal matrix with 
diagonal entries $t^{-w_0}, \ldots, t^{-w_n},$ for $t \neq 0.$  
Let $h \in \GL(n+1,K)$ be an upper 
triangular matrix.  If $g \in U$ then for any weight value $N$ we have 
\[
   \supp(h \cdot (\l gf)_{\geq N} ) = \supp((\l gf)_{\geq N}) 
\]
for almost all values of $t.$  
\end{lemma}

\begin{proof}
Note that the result of $h$ on a state  
$\x^{C_1}\wedge\cdots\wedge\x^{C_r}$ 
is a linear combination of states  
of the form $\x^{D_1}\wedge\cdots\wedge\x^{D_r},$ 
where $w\cdot D_i \geq w\cdot C_i.$  In particular every term other 
than $\x^{C_1}\wedge\cdots\wedge\x^{C_r}$ has weight strictly 
larger than $w \cdot (C_1 + \cdots + C_r).$  

Let $\tilde{f} = gf.$  
We compute
\begin{eqnarray*}
h \lambda \tilde{f} &=& h\lambda \sum_{M} \tilde{f}_{M} \\
 &=& h \sum_{M} t^{-M}\tilde{f}_{ M} \\
 &=& \sum_{M \geq  N} t^{-M} h \cdot \tilde{f}_{M} + 
     \sum_{M' < N} t^{-M'}h\cdot\tilde{f}_{M'}   \\
 &=& \left[ h \cdot (\l g f)_{\geq N} \right]  + 
	\left[ h \cdot (\l g f)_{< N} \right]
\end{eqnarray*}

Now suppose that $h \cdot ( (\l gf)_{\geq N} )$ 
contains a term 
$c(t) \x^{C_1}\wedge\cdots\wedge\x^{C_r}$ appearing with the 
coefficient $c(t) \neq 0,$ which we consider as a 
Laurent polynomial in the variable $t.$  
Let $c'(t)$ be the coefficient of 
the same term in the right-hand sum.
The degree of $c'(t)$ as a Laurent polynomial in $t$ is 
strictly greater than $-N,$ while that of $c(t)$ is 
at most $-N.$  Hence the two cannot cancel as polynomials
($c(t) + c'(t) \neq 0$) and the state 
$\x^{C_1}\wedge\cdots\wedge\x^{C_r}$ lies in the support 
of $h \cdot (\l g f)$ for infinitely many values of $t$ 
(recall that the ground field $K$ is infinite).  
By observations (2) and (3) above we see 
it lies in the support of $\l g f.$  This proves $\subseteq.$  

Conversely, if $c t^{-M} \x^{C_1}\wedge\cdots\wedge\x^{C_r}$ 
is a term of $(\l g f)_{\geq N}$ with weight 
$M = w \cdot (C_1 + \cdots + C_r) \geq N,$ where 
$c \in K \setminus \{0\},$ then 
its coefficient in $h \cdot (\l g f)_{\geq N}$ is a 
Laurent polynomial in $t$ which has $ct^{-N}$ as the only term 
with that power of $t$ occurring.  Hence it is not zero for 
infinitely many $t.$  There are only finitely many terms, so 
we get the other containment, $\supseteq.$  
\end{proof}

Now we are in a position to prove an infinitesimal version 
of Galligo's Theorem.  The idea is as follows.  We take an 
ideal truncated in a large degree, say $I = (f_1, \ldots, f_r).$  
What we want is to deform the ideal $I$ in generic coordinates 
to its initial ideal (given some monomial order), and show 
that as we get infinitesimally close to the initial ideal, we 
have something that is Borel-fixed.  Deforming to the initial 
ideal is done by acting by a diagonal matrix with diagonal 
entries $t^{-w_0}, \ldots, t^{-w_n},$ where 
$w = (w_0, \ldots, w_n)$ induces our monomial order.  Thus we 
will have a family of ideals parametrized by the variable $t,$ 
and taking the limit as $t \rightarrow 0$ gives the initial ideal.  

We see that as $t$ gets small, the monomials largest in the 
term order begin to dominate (they have $t$ coefficients with 
the smallest negative powers).  The monomials next highest 
in the monomial order will then govern the 
\emph{first order} behavior of the one parameter family.  
Specifically, if one takes the highest wedge product of the 
defining polynomials of a member of this family of ideals, 
then there is a unique term of highest weight (a fact we 
exploited to prove Galligo's theorem);  however there 
may be many terms with the next highest weight, and it 
is these terms that will dominate to first order.  The 
fact that there may be many terms of second highest weight 
presents a stumbling block.  When we were just interested 
in the unique term with highest weight we could argue that 
after acting by an upper triangular matrix we could not have 
produced a new term with higher weight, since we had already picked 
the largest one possible.  Now we need to controll the terms 
with the second highest weight.  However, 
though there may be many,   
\ref{support equality} at least gurantees that the set 
of these states remains invariant.  

The final problem we might encounter is that we don't really know 
what happens to the coefficients of the second highest weight 
terms after acting by an upper triangular matrix.  To 
remedy this we will use lemma 
\ref{almost Borel-fixed = Borel-fixed} which essentially 
says these coefficients are a red herring.  
Now on to the theorem.

\begin{theorem}\label{thm: ifinitesimal Galligo}
Let $z$ be a point on the Hilbert scheme corresponding 
to the subscheme $Z \subseteq \P^{n},$ and let 
$m$ be the Gotzmann number for the Hilbert polynomial 
of $Z.$  
Let $I$ be its defining ideal 
truncated at the degree $m.$  
Fix a weight vector 
$w = (w_0, \ldots, w_n)$ which distinguishes monomials 
in degrees at least up to $rm.$  As before set $\l = \l_w(t)$ to 
be the diagonal matrix with entries 
$t^{-w_0}, \ldots, t^{-w_n}.$  
Let $f_1, \ldots, f_r$ be a basis for $I$ (and 
thus a linear basis for $I_m$).  
Finally 
let $U = U(f),$ where 
$f = f_1\wedge\cdots\wedge f_r \in \wedge^r S_m$ 
(as defined above).
Then for $g \in U,$ the path on $\Hilb$ defined by the 
one-parameter family of ideals 
$\{ \l_w(t) gI\}$ has as limit as $t \rightarrow 0$ a 
Borel-fixed point, and the tangent vector to this 
path at that point is an eigenvector for the Borel 
group of upper triangular matrices.
\end{theorem}

\begin{proof}
\newcommand\tf{\tilde{f}}
We have $gI = (gf_1, \ldots, gf_r).$  
Let $\tf_1, \ldots, \tf_r$ be a new basis for 
$gI$ where 
$\init_w (\tf_i) = \x^{A_i}$ and this term appears 
in no other $\tf_j,$ so that the initial ideal of 
$gI$ is $(\x^{A_1}, \ldots, \x^{A_r}).$
We already know that this is Borel-fixed (Galligo's theorem). 
Let $\F$ be the filter of 
exponent vectors $\{A_i\}$ and $\R$ the 
order ideal of all other exponent vectors in degree $m.$  
Thus we can write
\[ 
   \tf_i = \x^{A_i} + \sum_{B \in \R} c_{A_i,B} \x^B
\] 
and therefore
\[ 
   \lambda \tf_i = t^{-w\cdot A_i}\x^{A_i} + 
   \sum_{B \in \R} c_{A_i,B} t^{-w\cdot B}\x^B 
\] 
For $t \neq 0$ then we find that $\lambda g I$ is 
generated by $(f_1', \ldots, f_r'),$ where for each 
$i$ we set 
\begin{equation}\label{eqn: f_i'} 
   f_i'\ :=\ t^{w\cdot A_i}\lambda \tf_i\ =\ \x^{A_i} + 
   \sum_{B \in \R} c_{A_i,B} t^{w\cdot(A_i-B)}\x^B. 
\end{equation} 
Among the set of all differences $A_i-B$ with 
$c_{A_i,B} \neq 0$ choose one 
$K = A_i-B$ such that $w\cdot K$ is minimal.  
Let 
\[
   \F' = \F \setminus \{ A_i \in \F \mid B=A_i+K \in \R 
     \mathrm{\ and\ } c_{A_i,B} \neq 0 \} 
\]
and
\[ 
   \F'' = \F \cup \{ B \in \R \mid A_i=B+K \in \F 
     \mathrm{\ and\ } c_{A_i,B} \neq 0 \} .
\]
As $t \rightarrow 0$ the smallest powers of $t$ dominate and 
we see the tangent vector is given by setting to 
zero all powers of $t$ greater than $w\cdot K.$  
Thus the tangent vector (as an ideal in 
$S[\e]$) is given by the basis 
\[ 
   \{\x^{A_i} \mid A_i \in \F'\} \cup 
   \{\x^{A_i}+\e c_{A_i,A_i+K}\x^{A_i+K} 
   \mid A_i \in \F\setminus\F'\}. 
\]
Note that this is an eigenvector for the maximal torus by 
\ref{the torus-fixed directions}, 
of type $(\F', K)$ (see section 
\ref{The Borel eigenvectors} for the 
definition of \emph{type}).    
By lemma \ref{almost Borel-fixed = Borel-fixed} what we 
need to show is that after acting by an upper triangular 
matrix we get a vector with the same type.  To do 
this return momentarily to the ideal 
$\l g I = (f_1', \ldots, f_r'),$  for $t \neq 0.$  
From equation \ref{eqn: f_i'} we see that after 
expanding $f_1'\wedge\cdots\wedge f_r'$ we will have 
$\x^{A_1}\wedge\cdots\wedge\x^{A_r}$ as the highest 
weight term, with weight $N = w\cdot(A_1+\cdots+A_r),$ 
and the second highest weight occurring is $N_1 = w\cdot(A+K).$  
Specifically, if 
$\a = \x^{A_1}\wedge\cdots\wedge\x^{A_r},$  
and if for 
$A_i \in \F\setminus\F'$ we set 
$\a_i = \x^{A_1}\wedge\cdots\wedge\x^{A_i+K}
\wedge\cdots\wedge\x^{A_r},$ 
(that is replace the monomial 
$\x^{A_i} \in \F\setminus\F'$ with 
$\x^{A_i+K}$), then we have 
\[ 
   f_1'\wedge\cdots\wedge f_r' = 
   \a + \sum_{A_i\in\F\setminus\F'} c_{A_i,A_i+K}t^{w\cdot K}\a_i + 
   (\mathrm{terms\ of\ lower\ weight}).
\]
Lemma \ref{support equality} gives us for any upper-triangular 
matrix $h$ that 
\[ 
   \supp(h \cdot(\lambda g f)_{\geq N_1}) = 
   \supp(\lambda gf)_{\geq N_1} 
\]
for almost all values of $t.$  Thus this holds for 
$t$ in a Zariski open subset of $\mathbf{A}^1\setminus 0.$  
Letting $t \rightarrow 0$ we see 
the tangent vector 
which has type $(\F', K)$ still has type 
$(\F', K)$ after acting by $h.$  Since $h$ was 
arbitrary, lemma \ref{almost Borel-fixed = Borel-fixed} 
says this vector is an eigenvector for the 
Borel subgroup of upper-triangular matrices. 
\end{proof}

Some comments are in order.  First we should note that the 
open subset $U$ in theorem 
\ref{thm: ifinitesimal Galligo} is smaller than that 
used for Galligo's theorem.  Thus ``generic coordinates'' 
has a stricter interpretation here.  
That said, we could have defined a larger open set on 
which the theorem still holds, but doing so drastically 
reduces a considerable degree of clarity. 

Second we should comment on the choice of 
weight vector $w.$  In theorem 
\ref{thm: ifinitesimal Galligo} we chose 
$w$ to distinguish monomials up to the 
large degree $rm,$ where $m$ is the Gotzmann number, 
and $r = \dim I_m.$  First off we could of chose 
the $m$ simply as the degree of definition of the 
initial ideal.  We simply chose to avoid over complicating 
the statement.  Second we choose the large 
degree $rm$ to 
ensure that states with distinct weights are weighted 
with distinct powers of the paremetrizing variable $t.$  
However any weight vector inducing the term order only 
up to degree $m$ already induces the term order, in the 
sense that the initial ideal with respect to the weight 
vector is the same as that with respect to the term order.  
Thus our condition on $w$ is considerably more strict.  
Put another way, given an ideal in generic coordinates, 
we can define the \emph{first order Gr\"{o}bner fan} by 
taking the open chambers to be those weight vectors 
producing the same Borel eigenvector.  This 
fan is finer than the typical Gr\"{o}bner fan.  Hence 
distinct weight vectors which induce the same term order 
may still produce different Borel eigenvectors.  
A weight vector that lies on a wall of the first order 
Gr\"{o}bner fan, but in an open chamber of the 
typical Gr\"{o}bner fan, wil not give a 
Borel eigenvector.

\bibliography{refs}

\begin{thebibliography}{Mum66}

\bibitem[Bay82]{bayer}
Dave Bayer.
\newblock The division algorithm and the {H}ilbert scheme.
\newblock {\em Thesis}, 1982.

\bibitem[BS87]{bayer/stillman}
David Bayer and Michael Stillman.
\newblock A theorem on refining division orders by the reverse lexicographic
  order.
\newblock {\em Duke Math. J.}, 55(2):321--328, 1987.

\bibitem[Eis95]{eisenbud}
David Eisenbud.
\newblock {\em Commutative algebra}, volume 150 of {\em Graduate Texts in
  Mathematics}.
\newblock Springer-Verlag, New York, 1995.
\newblock With a view toward algebraic geometry.

\bibitem[EK90]{eliahou/kervaire}
Shalom Eliahou and Michel Kervaire.
\newblock Minimal resolutions of some monomial ideals.
\newblock {\em J. Algebra}, 129(1):1--25, 1990.

\bibitem[Gal79]{galligo}
Andr{\'e} Galligo.
\newblock Th\'eor\`eme de division et stabilit\'e en g\'eom\'etrie analytique
  locale.
\newblock {\em Ann. Inst. Fourier (Grenoble)}, 29(2):vii, 107--184, 1979.

\bibitem[Got78]{gotzmann}
Gerd Gotzmann.
\newblock Eine {B}edingung f\"ur die {F}lachheit und das {H}ilbertpolynom eines
  graduierten {R}inges.
\newblock {\em Math. Z.}, 158(1):61--70, 1978.

\bibitem[Har77]{hartshorneAG}
Robin Hartshorne.
\newblock {\em Algebraic geometry}.
\newblock Springer-Verlag, New York, 1977.
\newblock Graduate Texts in Mathematics, No. 52.

\bibitem[HS04]{haiman/sturmfels}
Mark Haiman and Bernd Sturmfels.
\newblock Multigraded {H}ilbert schemes.
\newblock {\em J. Algebraic Geom.}, 13(4):725--769, 2004.

\bibitem[Mac27]{macaulay}
F.S. Macaulay.
\newblock Some properties of enumeration in the theory of modular systems.
\newblock {\em Proc. London Math. Soc.}, (26):531--555, 1927.

\bibitem[MR99]{marinari/ramella}
Maria~Grazia Marinari and Luciana Ramella.
\newblock Some properties of {B}orel ideals.
\newblock {\em J. Pure Appl. Algebra}, 139(1-3):183--200, 1999.
\newblock Effective methods in algebraic geometry (Saint-Malo, 1998).

\bibitem[Mum66]{mumford}
David Mumford.
\newblock {\em Lectures on curves on an algebraic surface}.
\newblock With a section by G. M. Bergman. Annals of Mathematics Studies, No.
  59. Princeton University Press, Princeton, N.J., 1966.

\bibitem[PS05]{peeva/stillman}
Irena Peeva and Mike Stillman.
\newblock Connectedness of {H}ilbert schemes.
\newblock {\em J. Algebraic Geom.}, 14(2):193--211, 2005.

\bibitem[Sne99]{snellman}
Jan Snellman.
\newblock On some partial orders associated to generic initial ideals.
\newblock {\em S\'em. Lothar. Combin.}, 43:Art.\ B43h, 23 pp. (electronic),
  1999.

\bibitem[Sta97]{stanley}
Richard~P. Stanley.
\newblock {\em Enumerative combinatorics. {V}ol. 1}, volume~49 of {\em
  Cambridge Studies in Advanced Mathematics}.
\newblock Cambridge University Press, Cambridge, 1997.
\newblock With a foreword by Gian-Carlo Rota, Corrected reprint of the 1986
  original.

\end{thebibliography}
\bibliographystyle{alpha}

\end{document}